\documentclass[a4paper,12pt]{amsart}
\usepackage{ifthen}
\usepackage{graphicx}
\usepackage{mathrsfs}

\numberwithin{equation}{section}
\nonstopmode
\newtheorem{thm}{Theorem}
\newtheorem{cor}[thm]{Corollary}
\newtheorem{lem}[thm]{Lemma}

\theoremstyle{definition}

\newtheorem{example}[thm]{Example}

\newenvironment{pf}[1][]{%
 \vskip 3mm
 \noindent
 \ifthenelse{\equal{#1}{}}%
  {{\slshape Proof. }}%
  {{\slshape #1.} }%
 }%
{\qed\bigskip}

\newcounter{alphabet}
\newcounter{tmp}

\newcommand{\an}{{\operatorname a}}

\newcommand{\C}{{\mathbb C}}
\newcommand{\D}{{\mathbb D}}

\newcommand{\M}{{\mathbf M}}

\newcommand{\N}{{\mathbf N}}

\newcommand{\R}{{\mathbb R}}

\newcommand{\T}{{\mathbb T}}

\newcommand{\bD}{{\overline{\mathbb D}}}

\renewcommand{\Re}{{\operatorname{Re}\,}}

\newcommand{\id}{{\operatorname{id}}}

\newcommand{\inv}{^{-1}}

\newcommand{\QC}{{\operatorname{QC}}}

\newcommand{\teich}{{\mathcal{T}}}

\newcommand{\aand}{{\quad\text{and}\quad}}
\newcommand{\zb}{{\bar{z}}}% z-bar

\newcommand{\esssup}{{\operatorname{esssup}\,}}

\newcounter{minutes}
\divide\time by 60
\newcounter{hours}
\multiply\time by 60
\addtocounter{minutes}{-\time}
\begin{document}
\bibliographystyle{amsplain}
\title{
A construction of trivial Beltrami coefficients
}

\def\thefootnote{}
%\footnotetext{
%\texttt{\tiny File:~\jobname .tex,
%          printed: \number\year-\number\month-\number\day,
%          \thehours.\ifnum\theminutes<10{0}\fi\theminutes}
%}
%\makeatletter\def\thefootnote{\@arabic\c@footnote}\makeatother

\author[T. Sugawa]{Toshiyuki Sugawa}
\address{Graduate School of Information Sciences,
Tohoku University, Aoba-ku, Sendai 980-8579, Japan}
\email{sugawa@math.is.tohoku.ac.jp}
\begin{abstract}
A measurable function $\mu$ on the unit disk $\mathbb{D}$ of the complex plane
with $\|\mu\|_\infty<1$ is sometimes called a Beltrami coefficient.
We say that $\mu$ is trivial if
it is the complex dilatation $f_\zb/f_z$ of a quasiconformal automorphism
$f$ of $\mathbb{D}$ satisfying the trivial boundary condition $f(z)=z,~|z|=1.$
Since it is not easy to solve the Beltrami equation explicitly,
to detect triviality of a given Beltrami coefficient is a hard problem,
in general.
In the present article, we offer a sufficient condition for a Beltrami coefficient
to be trivial.
Our proof is based on Betker's theorem on L\"owner chains.
\end{abstract}
%\thanks{
%The authors were supported in part by JSPS Grant-in-Aid for 
%Scientific Research (B) 22340025.
%}
\maketitle

\section{Introduction}
Let $k$ be a number with $0\le k<1.$
A homeomorphism $f$ of a plane domain $\Omega$ onto another $\Omega'$ is
called {\it $k$-quasiconformal} if $f$ has locally integrable distributional derivatives
$f_z=(f_x-if_y)/2$ and $f_\zb=(f_x+if_y)/2$ satisfying $|f_\zb|\le k|f_z|$
almost everywhere (a.e.) on $\Omega.$
If $k$ is not specified, $f$ is simply called quasiconformal.
It is well known that a homeomorphism $f$ is $k$-quasiconformal if and only if $f$ is ACL
(Absolutely Continuous on Lines), i.e., 
$f(x+iy)$ is absolutely continuous in $x$ for almost every (a.e.) $y$ and absolutely
continuous in $y$ for a.e.~$x$ for any rectangle $a\le x\le b,~ c\le y\le d$ contained in $\Omega$
and the partial derivatives satisfy $|f_\zb|\le k|f_z|$ a.e.~on $\Omega.$
See \cite{Ahlfors:qc} or \cite{LV:qc} for basic information about quasiconformal mappings.
It is known that $f_z\ne0$ a.e.~for a quasiconformal mapping $f.$
The ratio $\mu[f]=f_\zb/f_z$ is thus a well-defined measurable function on $\Omega$
with $\|\mu[f]\|_\infty\le k$ and called the {\it complex dilatation} of $f.$
Throughout the article, $\D$ will stand for
the open unit disk $\{z: |z|<1\}$ of the complex plane $\C.$ 
We denote by $\M(\D)$ the open unit ball of the Lebesgue space $L^\infty(\D).$
In other words, $\M(\D)$ consists of complex-valued measurable functions $\mu$ on $\D$ with
$\|\mu\|_\infty=\esssup_{z\in\D}|\mu(z)|<1.$
(Here and hereafter, we will always take a Borel measurable representative for
an element of the Lebesgue space.)
A member $\mu$ of $\M(\D)$ is often called a {\it Beltrami coefficient}.
Let $\QC(\D)$ stand for
the set of quasiconformal automorphisms $f$ of the unit disk $\D$
normalised so that $f(1)=1, f(i)=i$ and $f(-1)=-1.$
It follows from the measurable Riemann mapping theorem that there exists
a unique solution $f\in\QC(\D)$ to the Beltrami equation $f_z=\mu f_\zb$ a.e.~on $\D$
for a given $\mu\in\M(\D).$
We denote by $w^\mu$ the solution $f\in\QC(\D).$
A Beltrami coefficient $\mu\in\M(\D)$ is called {\it trivial} if $w^\mu=\id$ on 
the unit circle $\T=\partial\D.$
Set $\QC_0(\D)=\{f\in\QC(\D): f=\id~\text{on}~\partial\D\}$ and let
$\M_0(\D)$ the set of trivial Beltrami coefficients on $\D.$
Obviously, $\M_0(\D)$ corresponds to $\QC_0(\D)$ through the mapping $\mu\mapsto w^\mu.$
Note that $\QC(\D)$ is a group with composition as its group operation
and that $\QC_0(\D)$ is a normal subgroup of $\QC(\D).$
The {\it universal Teichm\"uller space} $\teich$ is defined as the quotient space
$\QC(\D)/\QC_0(\D).$
See \cite{Sergeev:ut} for basic information about the universal Teichm\"uller space.
It is also defined to be the quotient space $\M(\D)/\sim$ with respect to the
equivalence relation $\mu_1\sim\mu_2~\Leftrightarrow~w^{\mu_1}=w^{\mu_2}$ on $\T.$
Let $\pi:\M(\D)\to\teich$ be the natural projection.
By the Bers embedding, $\teich$ is realised as a bounded contractible domain in a complex
Banach space $B$ in such a way that $\pi:\M(\D)\to\teich$ is a holomorphic submersion
(see \cite{Nag:teich} for details).
Then the fiber $\pi\inv([0])$ over the basepoint $[0]=\pi(0)$ is nothing but $\M_0(\D).$
The tangent space $\N(\D)$ of $\M_0(\D)$ at $0$ is thus identified with
the kernel of the tangent map $d_{0}\pi:L^\infty(\D)\to B$ at $0.$
Each member of $\N(\D)\subset L^\infty(\D)$ is called {\it infinitesimally trivial}.
Teichm\"uller's lemma asserts that $\nu\in L^\infty(\D)$ is infinitesimally trivial if and only if
\begin{equation}\label{eq:it}
\iint_\D \nu(z)\varphi(z)dxdy=0\quad\text{for all}~\varphi\in A^1(\D),
\end{equation}
where $A^1(\D)$ denotes the Bergman space on $\D;$ namely,
the Banach space of integrable holomorphic functions on $\D.$
We refer to a monograph \cite{GL:teich} by Gardiner and Lakic
(more specifically, Theorem 6 in Chapter 6 therein)
for a proof and some background.
Earle and Eells \cite{EE67} proved that $\M_0(\D)$ is a contractible $C^0$-submanifold of $\M(\D).$
However, not much is known for the structure of $\M_0(\D)$ so far.
Difficulty in the study of $\M_0(\D)$ comes partly from the fact that there is no
convenient criterion for a Beltrami coefficient on $\D$ to be trivial.
It would be thus helpful to provide a method to construct a rich family
of trivial Beltrami coefficients.
One of such constructions is as follows.
Let $\Gamma$ be a torsion-free Fuchsian group acting on $\D$
such that the quotient Riemann surface $\D/\Gamma$ is biholomorphically
equivalent to the thrice-punctured sphere $X_0=\C\setminus\{0,1\}.$
Let $\Delta$ be a fundamental domain in $\D$ for $\Gamma$
whose boundary is of area zero, and consider the Beltrami coefficient
$$
\nu(z)=\sum_{\gamma\in\Gamma}\mu(\gamma(z))\frac{\overline{\gamma'(z)}}{\gamma'(z)},
\quad z\in\D,
$$
for $\mu\in\M(\D)$ with $\mu=0$ on $\D\setminus\Delta.$
Then $\nu\in\M_0(\D)$ because the thrice-punctured
sphere $X_0$ is rigid (i.e., it has no moduli).
The author learnt this from Katsuhiko Matsuzaki and
would like to thank him for it.
This is an interesting construction but the resulting $\nu$ is complicated.

The main purpose of the present article is to give a simple sufficient condition
for $\mu\in\M(\D)$ to be trivial.
To this end, we recall some basic facts about Hardy spaces of harmonic and
holomorphic functions.
We refer the reader to the standard textbook \cite{Duren:hp} on this matter.

Let $h^\infty$ denote the complex Banach space consisting of bounded 
complex-valued harmonic functions
on $\D$ with the supremum norm $\|\cdot\|_\infty.$
It is well known that each function $u\in h^\infty$ has a radial limit
$u^*(e^{i\theta})=\lim_{r\to1^-}u(re^{i\theta})$ for a.e.~$\theta\in\R.$
Conversely, through the Poisson integral,
each essentially bounded (complex-valued) function $U\in L^\infty(\T)$
on the unit circle $\T$ has a unique harmonic extension
$u$ to $\D$ in such a way that $U=u^*$ a.e.~on $\T.$
Note that $\|U\|_\infty=\|u\|_\infty.$
Therefore, the mapping $\tau: h^\infty\to L^\infty(\T)$
defined by $\tau(u)=u^*$ is an isometric isomorphism between complex Banach spaces.
Note that the set $H^\infty$ of bounded holomorphic functions on $\D$
is a closed subspace of $h^\infty.$
Hence, $L_\an^\infty(\T):=\tau(H^\infty)$ is a closed subspace of
$L^\infty(\T)$ consisting of boundary values of bounded holomorphic functions on $\D.$

For a given $\mu\in\M(\D),$
$$
U_t(\zeta)=\zeta^{-2}\mu(e^{-t}\zeta),\quad \zeta\in\T,
$$
is a (Borel measurable) function in $L^\infty(\T)$ with $\|U_t\|_\infty\le \|\mu\|_\infty=:k<1$
for a.e.~$t\in[0,+\infty).$
Set $\psi_t=\tau\inv(U_t)\in h^\infty.$
In other words, $\psi_t(z)$ is the Poisson integral of $U_t(e^{i\theta})$ for $t\ge0.$
Moreover, the function $\Psi(z,t)=\psi_t(z)$ on $\D\times[0,+\infty)$
is harmonic in $z$ for a.e.~$t$ and measurable in $t$ for every $z$
with $|\Psi(z,t)|\le k.$
Conversely, if such a function $\Psi(z,t)=\psi_t(z)$ is given, then
$\mu(z)=(z/|z|)^2\psi_{-\log|z|}^*(z/|z|)$ belongs to $\M(\D)$ and
satisfies $\|\mu\|_\infty\le k.$
Our main result is now stated as follows.

\begin{thm}\label{thm:main}
Let $\mu\in\M(\D)$ and set $U_t(\zeta)=\zeta^{-2}\mu(\zeta e^{-t}),~\zeta\in\T,$
for $t\in[0,+\infty).$
If $U_t\in L_\an^\infty(\T)$ for a.e.~$t\ge0,$ then $\mu$ is a trivial Beltrami coefficient.
\end{thm}

We denote by $\M_\an(\D)$ the set of those $\mu\in\M(\D)$ which
satisfy the assumption in the theorem.
The theorem means the inclusion relation $\M_\an(\D)\subset\M_0(\D).$
We remark that the proof given below yields $w^\mu(0)=0$ for $\mu\in\M_\an(\D).$
We note that $\M_\an(\D)$ is a linear slice of $\M(\D)$ and
that each member $\nu$ of $\M_\an(\D)$ certainly satisfies \eqref{eq:it};
in other words, $\M_\an(\D)\subset\N(\D).$

\begin{example}
Let $N$ be a positive integer.
We choose essentially bounded measurable functions $a_j(t)$ in $t\ge0$ 
and bounded holomorphic functions $\varphi_j(z)$ in $|z|<1$
for $j=1,\dots,N$ so that
$$
\sum_{j=1}^N\|a_j\|_\infty\|\varphi_j\|_\infty<1.
$$
Then,
$$
\mu(z)=\sum_{j=1}^N a_j(-\log|z|)
\left(\frac{z}{|z|}\right)^{2}\varphi_j^*\left(\frac{z}{|z|}\right),
\quad z\in\D,
$$
satisfies $\|\mu\|_\infty<1.$
Since
$$
U_t(\zeta)=\zeta^{-2}\mu(e^{-t}\zeta)=\sum_{j=1}^N a_j(t)\varphi_j^*(\zeta),
\quad |\zeta|=1,
$$
extends to the bounded holomorphic function $\sum a_j(t)\varphi_j(z)$ on $\D$ for a.e.~$t\ge0,$
Theorem \ref{thm:main} implies that $\mu\in\M_0(\D).$
\end{example}

For example, when $N=2,~ a_1(t)=(1+i)/10$ for $t>\log2,~
a_1(t)=(-2+i)/10$ for $0<t<\log2, \varphi_1(z)=(\sin z/z)^2,
a_2(t)=e^{-i\log t}/5,~ \varphi_2(z)=(z-2/3)/(1-2z/3),$ we compute
$f=w^\mu$ numerically in Figure 1.
We observe that the boundary values are close to those of the identity map.
This picture was created by Hirokazu Shimauchi based on the method given in 
his joint paper \cite{PS16} with M.~Porter.
We note that their normalisation conditions are $f(0)=1$ and $f(1)=1$ for a quasiconformal
automorphism of $\D.$
By the property $w^\mu(0)=0$ for $\mu\in\M_\an(\D),$
their method works without renormalisation.

\begin{figure}[htbp]
\begin{minipage}{0.4\hsize}
\begin{center}
\includegraphics[height=0.25\textheight]{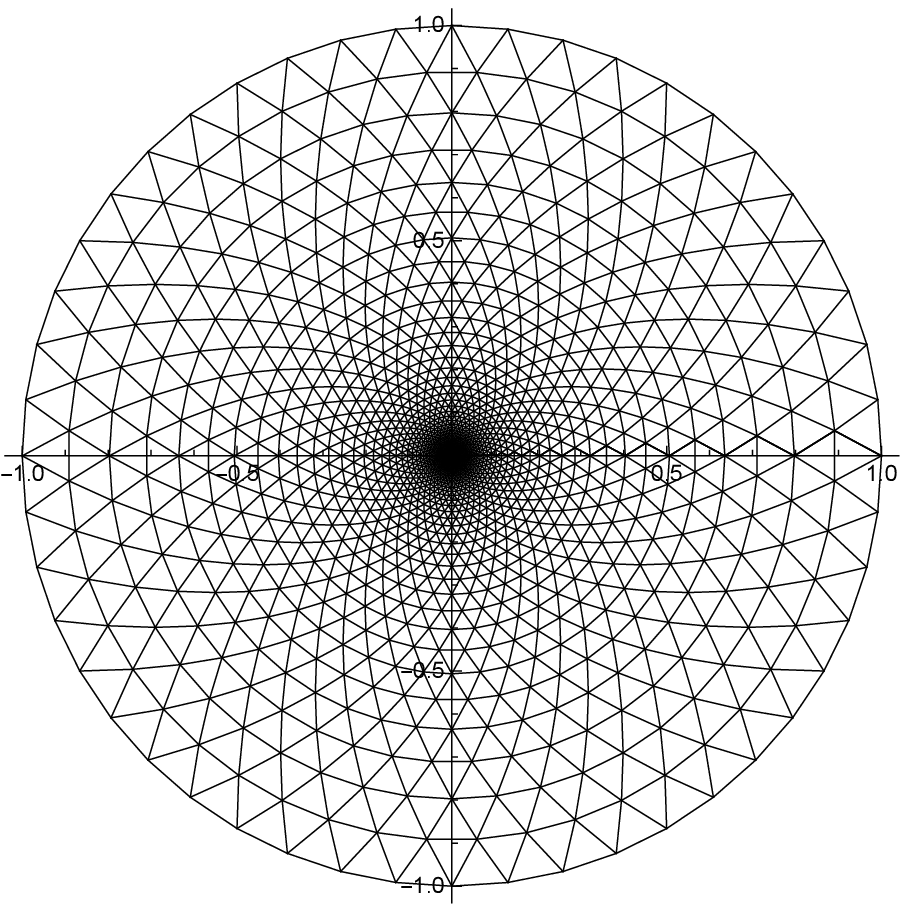}
\end{center}
\end{minipage}
\begin{minipage}{0.1\hsize}
\begin{center}
$\overset{f}{\longrightarrow}$
\end{center}
\end{minipage}
\begin{minipage}{0.4\hsize}
\begin{center}
\includegraphics[height=0.25\textheight]{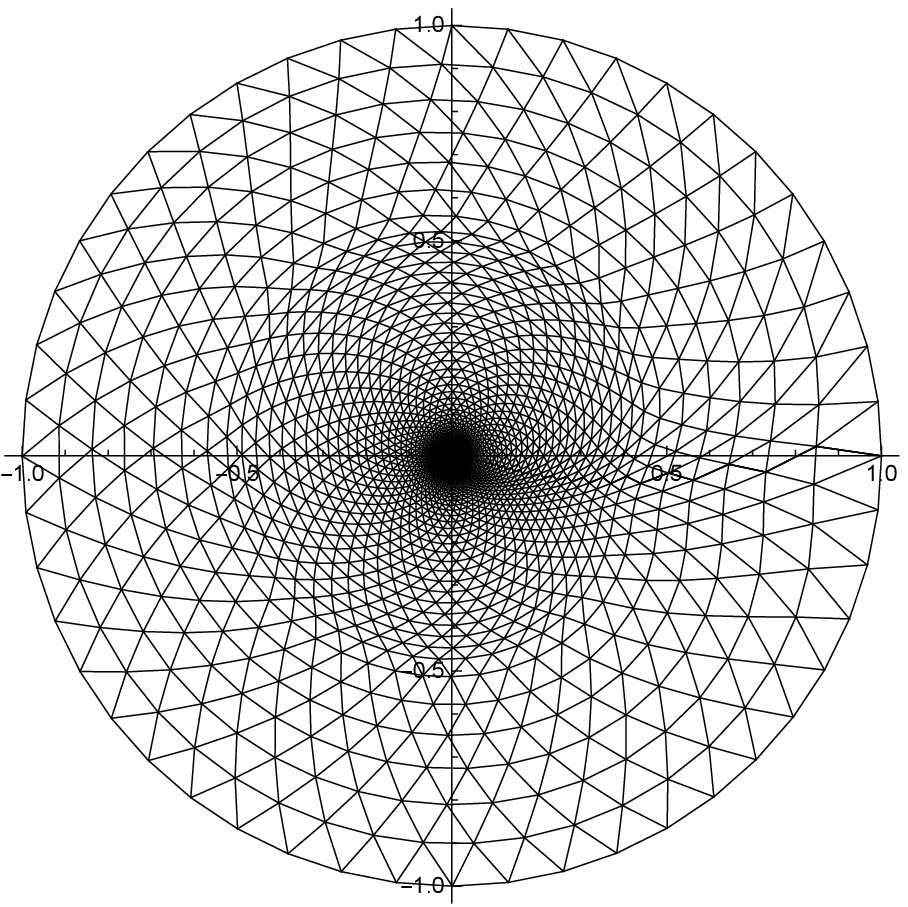}
\end{center}
\end{minipage}
\caption{A triangulation of the unit disk and its image under $f$}
\end{figure}

As a special case of Theorem \ref{thm:main}, we have the following result.

\begin{cor}
Let $\varphi$ be a bounded holomorphic function on $\D$ with $\|\varphi\|_\infty<1.$
Then $\mu(z)=z^2\varphi(z)$ is a trivial Beltrami coefficient on $\D.$
\end{cor}

\begin{pf}
Since $U_t(\zeta)=\zeta^{-2}\mu(\zeta e^{-t})=e^{-2t}\varphi(\zeta e^{-t})$
extends to the holomorphic function $e^{-2t}\varphi(ze^{-t})$ on $\D$ for each $t\ge0,$
the assertion follows from Theorem \ref{thm:main}.
\end{pf}

\section{Betker's theorem and proof of the main theorem}

In the proof of our main theorem, Betker's idea in \cite{Bet92} plays a cruicial role.
We now recall Betker's result.
Let $\omega_t(z)=\omega(z,t),~ t\ge0,$ be a family of holomorphic functions on
the unit disk $\D.$
It is called an {\it inverse L\"owner chain} if
\begin{enumerate}
\item[(i)]
$\omega(0,t)=w_0$ is independent of $t,$
\item[(ii)]
$b(t)=\omega_t'(0)$ is locally absolutely continuous in $t\ge0$
and $b(t)\to 0$ as $t\to+\infty,$
\item[(iii)]
$\omega_t:\D\to\C$ is univalent and $\omega_s(\D)\supset \omega_t(\D)$ whenever
$0\le s\le t.$
\end{enumerate}
We say that $q(z,t)=q_t(z)$ is a measurable family of holomorphic functions on $\D$
if it is holomorphic in $z\in\D$ for a.e.~$t\in[0,+\infty),$ and measurable in $t\ge0$
for each $z\in\D.$
As in the case of usual L\"owner chains, one can show that
$\omega(z,t)$ is absolutely continuous in $t$ for each $z$ and it
satisfies the partial differential equation
\begin{equation}\label{eq:il}
\dot \omega(z,t)=-z\omega'(z,t)q(z,t),\quad z\in\D,~ \text{a.e.}~t\ge0,
\end{equation}
for a measurable family of holomorphic functions $q(z,t)=q_t(z)$ on $\D$
with $\Re q_t>0$ on $\D$ for a.e.~$t\ge0,$
where $\dot\omega=\partial\omega/\partial t$ and $\omega'=\partial\omega/\partial z.$
Conversely, suppose that a measurable family $q_t(z)=q(z,t)$ of holomorphic
functions on $\D$ with $\Re q_t>0$ is given.
If furthermore
$$
\int_0^\infty \Re q(0,t)dt=+\infty,
$$
then there is an inverse L\"owner chain $\omega_t(z)=\omega(z,t)$ satisfying \eqref{eq:il}
and $\omega(z,0)=z$ for $z\in\D$ (see \cite[Lemma 1]{Bet92}).
Note that $\omega(0,t)=\omega(0,0)=0$ for $t\ge0$ in this case.
The proof of our main theorem will be based on the following lemma, which was used
by Betker to prove his main result in \cite{Bet92}.

\begin{lem}[$\text{Betker \cite[Lemma 2]{Bet92}}$]\label{lem:bet}
Let $0\le k<1$ be a constant.
Suppose that $q(z,t)=q_t(z)$ is a measurable family of holomorphic functions on $|z|<1$ satisfying
$$
\left|\frac{q(z,t)-1}{q(z,t)+1}\right|\le k
$$
for $z\in\D$ and a.e.~$t\ge0.$
Then there exists an inverse L\"owner chain
$\omega(z,t)=\omega_t(z)$ such that $\omega(z,0)=z,~z\in\D,$
\begin{equation}\label{eq:PDE}
\dot \omega(z,t)=-z\omega'(z,t)q(z,t),\quad z\in\D, \text{a.e.}~t\ge0,
\end{equation}
$\omega_t$ has a continuous and injective extension $\omega_t^*$
to $\overline\D$ for each $t\ge0.$
Moreover, the following mapping $f:\C\to\C$ is $k$-quasiconformal:
\begin{equation}\label{eq:f}
f(z)=\begin{cases}
0,& \quad z=0, \\
\omega^*(z/|z|,-\log|z|)=\omega_{-\log|z|}^*(z/|z|),
& \quad 0<|z|\le 1, \\
z,&\quad |z|>1.
\end{cases}
\end{equation}
\end{lem}

We will also make use of the following lemma (see Theorem 5.2 of Chapter IV in \cite{LV:qc}).

\begin{lem}\label{lem:LV}
Let $f_n,~n=1,2,\dots,$ be a sequence of $k$-quasiconformal mappings on a plane domain $\Omega.$
Suppose that $f_n$ converges to a quasiconformal mapping $f$ on $\Omega$ locally uniformly
and that the complex dilatation $\mu_n$ of $f_n$ converges to a measurable function $\mu_\infty$ a.e.~on $\Omega.$
Then the complex dilatation $\mu$ of $f$ coincides with $\mu_\infty.$
\end{lem}

We are now ready to show our main theorem.

\begin{pf}[Proof of Theorem \ref{thm:main}]
Let $\mu\in\M(\D)$ with $\|\mu\|_\infty\le k<1$
and set $U_t(\zeta)=\zeta^{-2}\mu(e^{-t}\zeta)$ for $\zeta\in\T$
and $t\ge0.$
Suppose that $U_t=\psi_t^*$ for $\psi_t\in H^\infty$ for a.e.~$t\ge0.$
We define a measurable family of holomorphic functions on $\D$ by
$$
q(z,t)=\frac{1+\psi_t(z)}{1-\psi_t(z)},\quad z\in\D, t\in[0,+\infty).
$$
Since $|\psi_t|\le k$ on $\D$ for a.e.~$t\ge0,$ we see that
$|(q(z,t)-1)/(q(z,t)+1)|=|\psi_t(z)|\le k$ for $z\in\D$ and a.e.~$t\ge0.$
Lemma \ref{lem:bet} now yields that there is an inverse L\"owner chain
$\omega_t, 0\le t<+\infty,$ satisfying $\omega_0=\id_\D,$ the equation
\eqref{eq:PDE} and that $\omega_t$ is continuous and injective on $\overline\D.$
Moreover, the function $f$ given by \eqref{eq:f} is a $k$-quasiconformal mapping
of $\C$ onto itself.
Note that $f$ maps $\D$ onto itself and that $f=\id$ on $\T.$
We also note that $f(0)=0$ by Lemma \ref{lem:bet}.
Hence, the complex dilatation $\mu[f]$ of $f|_\D$ is trivial; namely,
$\mu[f]\in\M_0(\D).$
We now find a form of the complex dilatation of $f(z)=\omega^*(\zeta,t)$ with $z=\zeta e^{-t},~ t\ge0.$
Since it is not clear that $\omega^*(\zeta, t)=\omega_t^*(\zeta)$ is absolutely continuous in $\zeta\in\T$
for a.e.~$t\ge0,$ we have to take a standard approximation procedure.
Let $r_n,~n=1,2,\dots,$ be a strictly increasing sequence of positive numbers converging to $1$
and set
$$
q_n(z,t)=q(r_nz,t) \aand
\omega_n(z,t)=r_n\inv\omega(r_nz,t),\quad z\in\bD,~ t\ge0.
$$
Then $\omega_n$ is a solution to the inverse L\"owner equation
$\dot\omega_n(z,t)=-z\omega_n'(z,t)q_n(z,t)$
for $|z|<1/r_n$ and $t\in[0,+\infty)\setminus E_n,$ where $E_n$ is a set of
linear measure zero.
Put $E=\cup_n E_n.$
Then $E$ is again of linear measure zero.
Let $f_n$ be the $k$-quasiconformal mapping constructed with $\omega_n$ in
Lemma \ref{lem:bet}.
By introducing the logarithmic coordinates $w=t+i\theta=-\log z,$
we consider the function 
$$
F_n(w)=f_n(e^{-w})=\omega_n(e^{-i\theta},t)
%=r_n\inv\omega(r_ne^{-i\theta},t).
$$
on $\Re w>\log r_n.$
Then, in view of \eqref{eq:PDE}, we obtain
\begin{align*}
\mu[F_n](w)&=\mu[f_n](e^{-w})e^{w-\bar w}=e^{2i\theta}\mu[f_n](e^{-t-i\theta}) \\
&=\frac{\partial_tF_n+i\partial_\theta F_n}{\partial_tF_n-i\partial_\theta F_n} \\
&=\frac{\dot\omega_n(e^{-i\theta},t)+e^{-i\theta}\omega_n'(e^{-i\theta},t)}%
{\dot\omega_n(e^{-i\theta},t)-e^{-i\theta}\omega_n'(e^{-i\theta},t)} \\
&=\frac{-q_n(e^{-i\theta},t)+1}{-q_n(e^{-i\theta},t)-1}=\psi_t(r_ne^{-i\theta})
\end{align*}
for $\theta\in\R$ and $t\in(0,+\infty)\setminus E.$
Therefore,
$$
\mu[f_n](\zeta e^{-t})=\zeta^2\psi_t(r_n\zeta)\quad \text{for}~\zeta\in\T,  t\in(0,+\infty)\setminus E.
$$
Letting $n\to\infty,$ we observe that for each $t\in(0,+\infty)\setminus E,$
$\mu[f_n](\zeta e^{-t})$ tends to $\zeta^2U_t(\zeta)=\mu(\zeta e^{-t})$ for a.e.~$\zeta\in\T.$
In other words, for each $r\in(0,1)$ with $-\log r\notin E,$
$g_n(re^{i\theta})$ tends to $0$ for a.e.~$\theta\in[0,2\pi],$
where $g_n(z)=|\mu[f_n](z)-\mu(z)|.$
Since the complex dilatations are essentially bounded by $1,$
Lebesgue's convergence theorem ensures that
$$
\int_0^{2\pi}g_n(re^{i\theta})d\theta\to0 \quad (n\to\infty)
$$
for a.e.~$r\in(0,1).$
Hence, by the convergence theorem again together with the Fubini theorem, we have
$$
\iint_\D g_n(z)dxdy=\int_0^1\int_0^{2\pi}g_n(re^{i\theta})rd\theta dr\to 0
\quad(n\to\infty).
$$
In particular, $g_n\to0$ in measure on $\D$ as $n\to\infty.$
By a theorem of F.~Riesz, we conclude that $g_{n_j}\to0$ a.e.~on $\D;$
in other words, $\mu[f_{n_j}]\to \mu$ a.e.~on $\D$ as $j\to\infty$
for a subsequence $n_j$ (see, for instance, \cite[Theorem 2.2.5]{Bog:mt1}).
Hence, Lemma \ref{lem:LV} implies that the complex dilatation of $f$ is equal to $\mu.$
We now conclude that $\mu=\mu[f]\in\M_0(\D).$
\end{pf}

{\bf Acknowledgements.}
The author would like to thank Hirokazu Shimauchi for checking the result numerically
through many examples and for producing pictures given in Figure 1.
He also thanks Katsuhiko Matsuzaki for suggestive conversations on the present matter.

\def\cprime{$'$} \def\cprime{$'$} \def\cprime{$'$}
\providecommand{\bysame}{\leavevmode\hbox to3em{\hrulefill}\thinspace}
\providecommand{\MR}{\relax\ifhmode\unskip\space\fi MR }
\providecommand{\MRhref}[2]{%
  \href{http://www.ams.org/mathscinet-getitem?mr=#1}{#2}
}
\providecommand{\href}[2]{#2}

%\bibliography{papers}

\begin{thebibliography}{10}

\bibitem{Ahlfors:qc}
L.~V. Ahlfors, \emph{Lectures on {Q}uasiconformal {M}appings}, van Nostrand,
  1966.

\bibitem{Bet92}
{Th}. Betker, \emph{L\"owner chains and quasiconformal extensions}, Complex
  Var. \textbf{20} (1992), 107--111.

\bibitem{Bog:mt1}
V.~I. Bogachev, \emph{Measure {T}herory {V}ol. {I}}, Springer Verlag, Berlin
  and Heiderberg, 2007.

\bibitem{Duren:hp}
P.~L. Duren, \emph{Theory of {$H^p$} {S}paces}, Academic Press, New York and
  London, 1970.

\bibitem{EE67}
C.~J. Earle and J.~Eells Jr., \emph{On the differential geometry of
  {T}eichm\"uller spaces}, J. Analyse Math. \textbf{19} (1967), 35--52.

\bibitem{GL:teich}
F.~P. Gardiner and N.~Lakic, \emph{Quasiconformal {T}eichm\"uller {T}heory},
  Amer. Math. Soc., 2000.

\bibitem{LV:qc}
O.~Lehto and K.~I. Virtanen, \emph{Quasiconformal {M}appings in the {P}lane,
  2nd {E}d.}, Springer-Verlag, 1973.

\bibitem{Nag:teich}
S.~Nag, \emph{The {C}omplex {A}nalytic {T}heory of {T}eichm\"uller {S}paces},
  Wiley, New York, 1988.

\bibitem{PS16}
R.~M. Porter and H.~Shimauchi, \emph{Numerical solution of the {B}eltrami
  equation via a purely linear system}, Constr. Approx. \textbf{43} (2016),
  371--407.

\bibitem{Sergeev:ut}
A.~Sergeev, \emph{Lectures on universal {T}eichm\"uller space}, EMS Series of
  Lectures in Mathematics, European Mathematical Society (EMS), Z\"urich, 2014.

\end{thebibliography}
\end{document}